\documentclass[12pt]{article}
\usepackage[centertags]{amsmath}
\usepackage{color,latexsym,amsfonts,amssymb,mathrsfs,stmaryrd}
\usepackage{breakcites}
\usepackage{graphicx}
\usepackage{flushend,cuted}
\usepackage{bm}
\usepackage{tabularx}
\usepackage{indentfirst}
\usepackage{xparse}
\usepackage{tikz}
\usepackage{mdwlist}
\usepackage{tkz-graph}
\GraphInit[vstyle = Shade]
\usetikzlibrary[intersections,
positioning,
petri,
backgrounds,
fit,
decorations.pathmorphing,
arrows,
arrows.meta,
bending,
calc,
intersections,
through,
backgrounds,
shapes.geometric,
quotes,
matrix,
trees,
shapes.symbols,
graphs,
math,
patterns,
external,
scopes,
matrix,
lindenmayersystems,
shapes.callouts,
shapes.misc,
angles,
shapes.arrows,
shadings]

\newcommand{\blue}[1]{\textcolor{blue}{#1}}

\allowdisplaybreaks
\newtheorem{thm}{Theorem}[section]
\newtheorem{prop}[thm]{Proposition}

\newtheorem{lem}[thm]{Lemma}
\newtheorem{rem}[thm]{Remark}
\newtheorem{counterex}[thm]{Counterexample}

\DeclareMathOperator*{\argmin}{arg\,min}
\def\b{{\beta}}
\def \bpi {\mbox{\boldmath$\pi$}}
\def\d{\delta}
\def\vt{\vartheta}
\def\e{\varepsilon}
\def\ex{\mathbb{E}}
\def \g {\mathscr{G}}
\def\iid{{i.i.d.}}
\def \l{{\lambda}}
\def\prob{\mathbb{P}}
\def \p {{\bf P}}
\def \q {{\bf Q}}
\def\re{{\mathbb{R}}}
\def\s{{\sigma}}

\def \snn {\sum_{j=2}^n \Xi_j}
\def \snnn {\sum_{k=3}^n \Xi_k}
\def \tv {{\rm TV}}
\def \supp {{\rm supp}}
\def\var{{\rm Var}}
\def \cv {{\cal V}}
\def \cw {{\cal W}}
\def \cz {{\cal Z}}

\def \A {{\mathscr{A}}}
\def \B {{\mathscr{B}}}
\def\D {\Delta}
\def \F {\mathscr{F}}
\def \G {{\Gamma}}
\def \H {{\mathscr{H}}}
\def \K {\mathscr{K}}
\def \L {\mathscr{L}}
\def\M {\mathscr{M}}
\def \P {{\mathscr{P}}}

\def\R{\mathbb{R}}
\def\Z{\mathbb{Z}}
\def\N{\mathbb{N}}

\newcommand{\qed}{{\hfill\hbox{\vrule width 5pt height 5pt depth 0pt}}}
\newcommand{\beas}{\begin{eqnarray*}}
\newcommand{\enas}{\end{eqnarray*}}
\newcommand{\bea}{\begin{eqnarray}}
\newcommand{\ena}{\end{eqnarray}}
\newcommand{\eq}{\begin{equation}}
\newcommand{\en}{\end{equation}}

\def \[ {\left[}
\def\]{\right]}
\def\({\left(}
\def\){\right)}

\def\ignore#1{}
\def\Ref#1{(\ref{#1})}

\begin{document}

\title{A large sample property in approximating the superposition of \iid\ point processes}

\author{ Tianshu Cong\footnote{School of Mathematics and Statistics,
The University of Melbourne,
VIC 3010, Australia, E-mail: tcong1@student.unimelb.edu.au. Work supported by a Research Training Program Scholarship and a Cross-Disciplinary PhD Scholarship in Mathematics and Statistics at the University of Melbourne.}, 
 Aihua Xia\footnote{School of Mathematics and Statistics,
The University of Melbourne,
VIC 3010, Australia, E-mail: xia@ms.unimelb.edu.au. Work supported in part by the Belz fund, Australian Research Council Grants Nos DP150101459 and DP190100613.}
 \ and \  Fuxi Zhang\footnote{School of Mathematical Sciences, Peking University, Beijing 100871, China, E-mail: zhangfxi@math.pku.edu.cn. 
 Work supported in part by NSF of China 11371040.
}
\\
}

 \maketitle

\begin{abstract}

One of the main differences between the central limit theorem and the Poisson law of small numbers is that the former possesses the large sample property (LSP), i.e., the error of normal approximation to the sum of $n$
independent identically distributed (\iid) random variables is a decreasing function of $n$. 
Since 1980's, considerable effort has been devoted to recovering the LSP for the law of small numbers in discrete random variable approximation.
In this paper, we aim to establish the LSP for the superposition of \iid\ point processes. 



\vskip12pt \noindent\textit{Running title\/}:
Superposition of Point Processes
\end{abstract}

\section{Introduction}
\setcounter{equation}{0}

The central limit theorem states that the distribution of the sum $S : = \sum_{i=1}^nX_i$ of independent copies of a random variable $X$ with finite second moment, after being normalized, converges weakly to the standard normal distribution. The Berry-Esseen bound ensures that, if $X$ has the finite third moment, the error of the normal approximation, measured in the Kolmogorov metric, is not worse than $c/\sqrt{n}$, where $c$ is a constant determined by the distribution of $X$. In other words, the central limit theorem has the large sample property (LSP), i.e., the quality of the approximation improves as the sample size becomes large. The LSP can also be established for the functional central limit theorem measured in the L\'evy-Prokhorov distance \cite{BS80,H84,Ku85, F94,Utev86}. Moreover, Stein's method can be used to estimate the errors of diffusion approximation \cite{Barbour90}.

The Poisson law of small numbers, on the other hand, does not possess the LSP. More precisely, if $X_i$'s are independent indicator random variables with $\prob(X_i=1)=1-\prob(X_i=0)=p_i$ for each $i$, then the total variation distance between the distribution of $W=\sum_{i=1}^nI_i$ and the Poisson distribution with mean $\lambda:=\sum_{i=1}^np_i$ is of the order ${\Omega}\left(\lambda^{-1}\sum_{i=1}^n p_i^2\right)$ \cite{BH84}. In particular, if $p_i=p$ for all $i$, one can see that the quality of approximation does not improve when $n$ becomes large. This is due to the fact that a Poisson distribution has only one parameter while a normal distribution has two parameters. To recover the LSP, one has to introduce more parameters into the approximating distributions, e.g.,  signed compound Poisson measures, translated  Poisson, compound Poisson, negative binomial and polynomial birth-death distributions
\cite{Presman,Kru,Ce,BX99,BC04,Roelin07,BCL92,BP99,BX01}.

If we consider point processes rather than {nonnegative integer-valued} random variables, the counterpart is the superposition $\cv_n = \Xi_1 + \cdots + \Xi_n$ of point processes $\{\Xi_i:\ 1\le i\le n\}$. The pioneering work of Grigelionis \cite{Grigelionis}  demonstrates that the distribution of {the} superposition of independent sparse point processes on the carrier space $\R_+$ converges weakly to a Poisson process distribution. The same phenomenon can be established for the superposition of dependent sparse point processes on a general carrier space \cite{Goldman,Jagers72,Brown78,Kallenberg83}. The accuracy of Poisson point process approximation has been of considerable interest since 1970's \cite{Ser75,Brown78}. Stein's method for Poisson process approximation was subsequently established by \cite{Barbour88,BB92} for estimating the approximation errors and the method was further refined by \cite{BWX,Xia05b,CX04}. In the context of the aforementioned superposition of \iid\ point processes, no error estimates were studied until the last decade \cite{Schuhmacher,CX11} and these studies show that the Poisson point process approximation to the superposition of \iid\ point processes does not possess the LSP either. The aim of this note is to show that, by introducing more parameters into the approximating point process distribution, it is possible to recover a LSP in approximating the superposition of \iid\ point processes.

Given that a Poisson point process on a compact metric space can be viewed as a Poisson number $Z$ of \iid\ points in the space, a natural step of introducing more parameters into the approximating point process is to replace the Poisson number $Z$ by a random variable $N$ whose distribution is controlled by two or more parameters, such as the translated  Poisson \cite{BC04,Roelin07}, negative binomial \cite{BP99} and polynomial birth-death distributions \cite{BX01}. The family of approximating distributions we will consider in this note is the polynomial birth-death process distributions introduced in \cite{XZ08}. To quantify the difference between two point processes, as in \cite{Schuhmacher,CX11}, we use the Wasserstein distance $d_2$ initiated in \cite{BB92}. The formal statement of the main result is given in Theorem~\ref{mainresult}. Several applications are provided in Section~\ref{Section.Example} to illustrate the order of convergence in the LSP. Section~\ref{Section.Proof} is devoted to the proof of the main result.

\section{Preliminaries and the main result} \label{Section.Main}

{\em 1. Point processes.} For the reader's convenience, in this part, we collect some basic concepts and facts, and introduce a {\it partitional total variation distance} 
for comparing point processes under a partition of the carrier space. The basic concepts needed for this note are point process, reduced palm process \cite[Chapter~10]{Kallenberg83}, the Wasserstein distance $d_2$ \cite{BB92} and partition \cite{XZ12}.

Let $\G$ be a compact metric space with metric $d_0$ bounded by 1. Let $\B (\G)$ be the Borel $\s$-algebra induced by $d_0$. A {\em configuration} $\xi$ on $\G$ is a collection of finitely many particles located in $\G$. Equivalently, it can be represented as a non-negative integer-valued finite measure on $\G$. Denote by $|\rho|$ the total mass of a measure $\rho$. Therefore, we can write $\xi$ as $\sum_{i=1}^{|\xi|} \d_{x_i}$, where $\d_x$ is the Dirac measure at $x$. Let $\H$ be the set of all configurations on $\G$, and $\B (\H)$ be the $\s$-algebra generated by the mappings $\xi \mapsto \xi (C)$, $C \in \B (\G)$ \cite[p.~12]{Kallenberg83}. A  {\em point process} is a measurable mapping from a probability space into $(\H,\B (\H))$. We use $\xi, \eta, \cdots$ to stand for configurations, $\Xi$, $\cv, \cw, \cz, \cdots$ to stand for point processes, and $\p, \q, \L(\Xi), \L(\cv), \cdots$ to stand for the laws of point processes.

Let $\Xi$ be a point process with finite mean measure $\l (d x) := \ex \Xi (d x)$. The family of point processes $\{ \Xi_x : x \in \G \}$ are said to be the {\it reduced Palm processes} associated with $\Xi$ if for any measurable function $f : \G \times \H \rightarrow \re_+:=[0,\infty)$,
 $$
\ex \[ \int_\G f(x,\Xi-\d_x)\Xi(dx) \] = \int_\G  \ex f(x,\Xi_x)  \l(dx) ,
 $$
\cite[Chapter~10]{Kallenberg83}. Furthermore, suppose $\l^{[2]}(dx,dy) : = \ex \Xi(dx)(\Xi-\d_x)(dy)$ is finite, then one can define the {\it second order reduced Palm processes} $\{ \Xi_{xy}: x,y \in \G \}$ associated with $\Xi$ by
 $$
\ex \[ \iint_{\G^2}
f(x,y;\Xi-\d_x-\d_y)\Xi(dx)(\Xi-\d_x)(dy) \] = \iint_{\G^2} \ex
f(x,y;\Xi_{xy}) \l^{[2]} (dx,dy) ,
 $$
for any measurable function $f : \G^2 \times \H \rightarrow [0,\infty)$ \cite[Chapter~12]{Kallenberg83}.

\cite{BB92} introduce a Wasserstein distance $d_2$ for quantifying the difference between two probability measures $\p$, $\q$ on $(\H,\B (\H))$.  The metric is defined in two stages. First, for two finite measures $\rho_1$ and $\rho_2$ on $\G$, define  $$
d_1(\rho_1,\rho_2) :=\left\{\begin{array}{ll}
1,& \mbox{for }|\rho_1| \neq |\rho_2|,\\
0,& \mbox{for }|\rho_1|=|\rho_2|=0,\\
\sup_{u \in \K} \left| \bar{\rho}_1 (u) - \bar{\rho}_2 (u) \right| ,& \mbox{for }|\rho_1|= |\rho_2|>0,
\end{array}\right.
 $$
where $\bar {\rho} : = \rho / |\rho|$ is the normalized measure of $\rho$, $\K: = \{ u : |u(x)-u(y)|\le d_0(x,y) , \forall x,y\in\G \}$ and $\rho (u) : = \int_{\G} u d \rho$. In particular, by the Kantorovich-Rubinstein duality theorem \cite[Theorem 8.1.1]{Rachev}, for two probability measures  $\mu$ and $\nu$ on $\G$, $d_1 (\mu, \nu) = \inf_{X \sim \mu, Y \sim \nu} \ex d_0 (X,Y) $, where $X,\ Y$ are $\G$-valued $\B (\G)$-measurable random elements. For two configurations $\xi_1:=\sum_{i=1}^{|\xi_1|} \d_{x_{1i}},\xi_2:=\sum_{i=1}^{|\xi_2|} \d_{x_{2i}}\in \H$, {we have the duality representation $d_1(\xi_1,\xi_2)=\min_{\boldsymbol{\pi}}\left\lbrace n^{-1} \sum_{i=1}^nd_0
	(x_{1i},x_{2\boldsymbol{\pi}(i)})\right\rbrace  $} when $|\xi_1|=|\xi_2|=n$ and $1$ otherwise, where $\min$ is taken over all permutations of $\{1,2,\cdots,n\}$.  The metric $d_2$ is defined as
 $$
d_2 (\p, \q) : = \sup_{f \in \F } |\p (f) - \q (f)|
= \inf_{\cw \sim {\bf P}, \cz \sim {\bf Q}} \ex d_1 (\cw, \cz) ,
 $$
where $ \F: = \{ f: |f(\xi) - f(\eta)| \le d_1 (\xi, \eta), \ \forall \
\xi, \eta \in \H \} $, \label{pageofcapitalF}and the last equality is due to the duality theorem \cite[Theorem 8.1.1]{Rachev}.

For a partition $\g = \{ G_i:\ i\in I \} \subset \B (\G)$ of $\G$, where $I\subset\N:=\{1,2,3,\dots\}$ is a finite set, let $t_i\in\argmin_{x}\sup_{s \in G_i} d_0 (s, x)$, that is, $t_i \in \G$ is a point such that $d_0 (G_i,t_i) : = \sup_{s \in G_i} d_0 (s, t_i)$ is as small as possible, $i\in I$. We call $t_i$ a {\em center} of $G_i$. Let $d_0 (\g) := \max_{i\in I} d_0 (G_i,t_i)$. We call $\g$ an {\em $\e$-partition} of $\G$ if $d_0 (\g) \le \e$. Denote all $\e$-partitions of $\G$ by $\P_\e$.

For any partition $\g=\{ G_i:\ i\in I \} $, we define an {\it assembling mapping} $\M_\g$ as
 $$
\M_\g \circ \eta := \sum_{i\in I} \eta (G_i) \d_{t_i}.
 $$
The assembling mapping, when applied to a configuration $\eta$, shifts all particles of $\eta$ in $G_i\in\g$ to its center $t_i$.
For a point process $\cw$, we define the {\it partitional total variation distance} as 
 \eq \label{Eq.tv}
\tv_\g (\cw) :=  \max_{i\in I} d_{tv} (\L (\M_\g \circ \cw + \d_{t_i}) ; \L (\M_\g \circ \cw) ) ,
 \en
 where for two probability measures $\p$ and $\q$ on $\H$, $d_{tv}(\p,\q)=\sup_{A\in\B(\H)}|\p(A)-\q(A)|.$ We {write}
 $$
\vt_\e (\cw) := \inf_{ \g \in \P_\e}  \tv_\g (\cw) .
 $$

{\em 2. Polynomial birth-death point process.} As mentioned in the Introduction, there are various ways to introduce more parameters into the approximating point process for better accuracy of approximation. In this part, we collect the facts around the polynomial birth-death point process established in \cite{XZ08}.

For $a > 0$, $0 \le b < 1$, $\b \ge 0$, we define the polynomial birth-death distribution introduced in \cite{BX01} as
$$\pi_{a,b;\b}(i+1):={\cal C}\prod_{j=0}^i\frac{a+bj}{(j+1)(1+\b j)},\ i\in \Z_+:=\{0,1,2,\dots\},$$
where $${\cal C}:={\cal C}(a,b;\b)=\left(1+\sum_{i=0}^\infty\prod_{j=0}^i\frac{a+bj}{(j+1)(1+\b j)}\right)^{-1}.$$
The distribution can be viewed as the equilibrium distribution of the birth-death process with birth rates $\{a+bk:\ k\in\Z_+\}$ and death rates $\{k(1+\b (k-1)):\ k\in\N\}$. The
polynomial birth-death point process is given by
$$\cz := \sum_{n=1}^Z \d_{U_n},$$
where $Z; U_1, U_2, \cdots$ are independent, $Z \sim \pi_{a,b;\b}$, $U_n \sim \mu$ with $\mu$ being a probability measure on $\G$, $\forall n \ge 1$.
Denote $\L (\cz)$ by $\bpi_{a,b;\b;\mu}$.

{\em 3. Main Result.} Suppose $\Xi_1, \cdots, \Xi_n$ are independent and identically distributed point processes. In each of the following two cases, we can give the polynomial birth-death point process approximation of the superposition $\cv_n= \Xi_1 + \cdots + \Xi_n$ under the Wasserstein distance $d_2$. Denote
 $$
\l (dx) : = \ex \Xi_1 (d x), \ \ \l^{[2]} (dx, dy) : = \ex \Xi_1 (dx) (\Xi_1 - \d_x) (dy),
 $$
 $$
\theta_r := \ex \[ |\Xi_1| (|\Xi_1| - 1 ) \cdots (|\Xi_1| - r+1) \] , \ \ \ \forall r \in\N .
 $$

\noindent {\bf Case 1}. If $\var (|\Xi_1|) \ge \ex |\Xi_1|$, we take $b = \frac {\theta_2 - \theta_1^2}{ \theta_2 - \theta_1^2 + \theta_1}$ and $\b = 0$.

\noindent {\bf Case 2}. Either $\frac 1 2 \ex |\Xi_1| < \var (|\Xi_1|) < \ex |\Xi_1|$ or $\frac 1 2 \ex |\Xi_1| = \var (|\Xi_1|)$ and $\theta_3 > \theta_2 (\theta_1 - 1)$, we set $b = 0$ and
 $$
\b = \frac  {\theta_1^2 - \theta_2} {(n-1)\theta_1 \[  \theta_1 - 2 ( \theta_1^2- \theta_2) \]  + (\theta_3 + \theta_2 - \theta_1 \theta_2) }.
 $$

\begin{rem} Note that $\theta_2 - \theta_1^2 = \var \( |\Xi_1| \) - \ex |\Xi_1| $ and 
$$ \theta_1 - 2 ( \theta_1^2- \theta_2)  = 2 \[ \var (|\Xi_1|) - \frac 1 2 \ex |\Xi_1| \],$$ we have $0 \le b < 1$ and $\b > 0$ for $n>1+\frac{\theta_1\theta_2-\theta_2-\theta_3}{\theta_1(\theta_1-2(\theta_1^2-\theta_2))}\vee 0$. 
\end{rem}

In all cases, let
 \bea
a & = & n \[ (1-b)  \theta_1 +  \b \theta_2 + \b (n-1) \theta_1^2 \] , \label{Eq.a} \\
\mu(dx) & = & \frac {1 + \b \ex |\Xi_{1,x}| + \b (n-1) \theta_1 } {\theta_1 + \b \[ \theta_2 + (n-1) \theta_1^2 \]} \l (dx) , \label{Eq.mu}
 \ena
where $\Xi_{1,x}$ is the reduced Palm distribution of $\Xi_1$ at $x$. Our main result is as follows.

 \begin{thm} \label{mainresult}
For both cases above, there exists a constant $C$, depending on $\L (\Xi_1)$, and $0 < \rho < 1$ such that for any $\e>0$,
 \begin{eqnarray}
&&d_2 \( \L \( \sum_{k=1}^n \Xi_k \) , \bpi_{a,b ; \b; \mu} \) \nonumber\\
&&\le 2 \e + C \vt_\e \( \sum_{k=1}^{n-2} \Xi_k \) + Cn \prob \( \sum_{k=1}^{n-2} | \Xi_k | \le \rho (n-2) \theta_1 \)
\label{mainresultineq1}\\
&&\le 2 \e + C \vt_\e \( \sum_{k=1}^{n-2} \Xi_k \) + O\left(n^{-1/2}\right)\label{mainresultineq1.1}
 \end{eqnarray}
 for $n>2$ in Case 1 and $n>1+\frac{\theta_1\theta_2-\theta_2-\theta_3}{\theta_1(\theta_1-2(\theta_1^2-\theta_2))}\vee 0$ in Case 2.
In particular,
 \begin{eqnarray}
&&d_2 \( \L \( \sum_{k=1}^n \Xi_k \) , \bpi_{a,b ; \b; \mu} \)\nonumber\\
 &&\le  (C+1) 2 \e_{n-2} + C n \prob \( \sum_{k=1}^{n-2} | \Xi_k | \le \rho (n-2) \theta_1 \) \label{mainresultineq2}\\
&&\le  (C+1) 2 \e_{n-2} + O\left(n^{-1/2}\right),\label{mainresultineq2.1}
 \end{eqnarray}
valid for the same range of $n$ specified above, where $\e_n$  is defined as
 $$
\e_n = \inf\left\{v:\ 2v\ge \vt_v  \( \sum_{k=1}^n \Xi_k \)\right\}.$$
 \end{thm}

\begin{rem}{\rm The last terms of \Ref{mainresultineq1} and \Ref{mainresultineq2} are typically of order $O(e^{-cn})$ for a positive constant $c$ depending on the distribution of $\Xi_1$ \cite[Theorem~2.7]{CL06} but the remaining terms of \Ref{mainresultineq1} and \Ref{mainresultineq2} are typically of order no better than $O\left(n^{-1/2}\right)$.}
\end{rem}

 \begin{prop} \label{Prop.epsilon}
$\e_n$ is decreasing in $n$.
 \end{prop}
{\em Proof.} First, $\vt_\e (\sum_{i=1}^{n+1} \Xi_i) \le \vt_\e (\sum_{i=1}^n \Xi_i)$ since $TV_\g (\sum_{i=1}^{n+1} \Xi_i) \le TV_\g (\sum_{i=1}^n \Xi_i)$ for any partition $\g$. It follows that $\vt_{\e_n} (\sum_{i=1}^{n+1} \Xi_i) \le \vt_{\e_n} (\sum_{i=1}^n \Xi_i) \le 2 \e_n$. Noting that $\vt_\e ( \sum_{i=1}^{n+1} \Xi_i ) > 2 \e$ when $\e < \e_{n+1}$ and $\vt_\e (\sum_{i=1}^{n+1} \Xi_i) \le 2 \e$ when $\e > \e_{n+1}$, we obtain $\e_n \ge \e_{n+1}$. \qed

\begin{rem} {\rm Case 1 is known as {\it over-dispersion} \cite{Fad94}. It is shown in \cite{BHX98} that over-dispersion in statistics arising from natural phenomena is much more common than under-dispersion, i.e., $\var (|\Xi_1|) < \ex |\Xi_1|$\ignore{, ``coinciding with folk-lore and practice"}.}
\end{rem}

 \begin{rem}\label{meanchange} {\rm Let $\nu (dx)= \frac 1 \theta_1\l (dx)$ be the normalized distribution of $\l (dx)$.
Since $\mu(dx)$ is the normalized distribution of $\l (dx) + \frac { \b \ex |\Xi_{1,x}| }{1 + \b (n-1) \theta_1} \l (dx)$, we have
$$d_2 \( \bpi_{a,b ; \b; \mu} , \bpi_{a,b ; \b; \nu} \) \le d_1 (\mu, \nu) = O (n^{-1}).$$ Thus $\mu$ in Theorem~\ref{mainresult} can be replaced by $\nu$ at the cost of $O(n^{-1})$ being added to the upper bound.}
 \end{rem}

\ignore{ \begin{rem}
If the point process $\Xi_1$ is near space homogeneous, one can take $\e = O (n^{-r})$. Then $k$ has order $n^{rd}$, where $d$ is the dimension of $\G$. { Consequently, $\vt_\e (\Xi_1) = O (n^{rd/2})$. Why?} In general, if we set $u_1=\frac12\wedge\left(1-d_{tv} ( \L ( \M_\g \circ \Xi_1 + \vt_{t_i} ), \L (\M_\g \circ \Xi_1))\right)$, then arguing as in the proof of Proposition~4.6 of \cite{BX99} \blue{(cf. \cite{MR07})}, we have
 \blue{ $$
d_{tv} (\L (\M_\g \circ \Xi + \vt_{t_i}), \L (\M_\g \circ \Xi)) \le \frac 1 {\sqrt{nu_1}}.
 $$ }
Thus has order $n^{-1/2+ rd/2}$. By taking $r = \frac 1 {d+2}$, one has $\e_n = O (n^{- \frac 1 {d+2}} )$.
 \end{rem}}

\section{Examples} \label{Section.Example}

In this section, we demonstrate the use of Theorem~\ref{mainresult} in five applications: Bernoulli process, Bernoulli process with shifts, compound Poisson process, renewal process and entrances and exits of Markov process. For simplicity, except in subsection~\ref{sectcompoundPoisson}, we only consider point processes on the carrier space $\G=[0,1]$ with $d_0(x,y)=|x-y|$. Extension to any compact carrier space is a straightforward exercise. 

\subsection{Bernoulli process}\label{sectBernoulli}

As a warming up example, we consider a simple Bernoulli process $\Xi_1 = \sum_{i=1}^mI_i \d_{t_i}$, where $\{I_i:\ 1\le i\le m\}$ are independent Bernoulli random variables with $\prob(I_i=1)=1-\prob(I_i=0)=p_i\in (0,1)$, $\{t_i:\ 1\le i\le m\}\subset\G$ and $m$ is a finite positive integer. This is a typical case where the actual support space of the point process is a subset of the carrier space and it reminds us that the partition technique should not be applied blindly. The following theorem is a generalisation of \cite{XZ08} with the same order of convergence as that for the special case in \cite{XZ08}.

\begin{thm} \label{3.1} For \iid\  Bernoulli processes  $\{\Xi_k\}_{k \in \mathbb{N}}$, if $\sum_{i=1}^mp_i(1-p_i)>\frac12\sum_{i=1}^mp_i$, let $a,b;\b$ and $\mu$ be as defined in Case 2, we have

	$$
	d_2 \( \L \( \sum_{k=1}^n \Xi_k \) , \bpi_{a,b ; \b; \mu} \)=O\left(n^{-1/2}\right)\mbox{ as }n\to\infty.
	$$
\end{thm}

\begin{rem}\label{Bernoullinorole}
{\rm The distances amongst $\{t_i:\ 1\le i\le m\}$ play no role in the speed of convergence. }
\end{rem}

\noindent{\it Proof of Theorem~\ref{3.1}.} The support of
$\Xi_1$ is a reduced carrier space $\G_r:=\{t_i:\ 1\le i\le m\}$, so it suffices to consider the reduced carrier space $\G_r$ with partition $\g = \{\{t_i\}:\ 1\le i\le m\}$ and $\e=0$.
With the partition $\g$, $\Xi_1$ corresponds to the vector $\vec X = (I_1, \cdots , I_m)$. Let $\vec Y$ be the sum of $n-2$ independent copies of $\vec X$, and $e_j$ be the vector with value $1$ at the $j$-th component and $0$ otherwise. Then, by the independence,
\begin{equation}
d_{tv} (\L (\vec Y+ e_j), \L (\vec Y )) = d_{tv} (\L (Y_j + 1), \L (Y_j)).\label{bernoulli1}
\end{equation}
Noting that $Y_j$ is the sum of independent Bernoulli$(p_j)$ random variables, we have $Y_j\sim{\rm Binomial}(n-2,p_j)$, which implies
$$
d_{tv} (\L (Y_j + 1), \L (Y_j)) =\max_{0\le l\le n-2}{n-2\choose l}p_j^l (1-p_j)^{n-l} \le\frac{C}{\sqrt{(n-2) p_j(1-p_j)}}.
$$
Hence, it follows from \Ref{bernoulli1} that the second term of \Ref{mainresultineq1.1} is bounded by $O\left(n^{-1/2}\right)$.
\qed

\subsection{Bernoulli process with shifts}\label{sectBernoullishifts}

The aim of this example is to show that we may use a marked point process to get better approximation bounds.

Similar to the previous subsection, we define $\Xi_1 = \sum_{i=1}^mI_i \d_{\zeta_i}$, where $\{(I_i,\zeta_i)\}$ are independent with $I_i$ having Bernoulli distribution with $\prob(I_i=1)=1-\prob(I_i=0)=p_i\in (0,1)$, $\zeta_i$ taking values in $\G$, and $m$ is a fixed positive integer. 
	
\begin{thm} \label{3.2} If $\sum_{i=1}^mp_i(1-p_i)>\frac12\sum_{i=1}^mp_i$, let $a,b;\b$ and $\mu$ be as defined in Case 2, then for $\iid$ Bernoulli processes with shifts $\{\Xi_k\}_{k \in \mathbb{N}}$, 
	$$
	d_2 \( \L \( \sum_{k=1}^n \Xi_k \) , \bpi_{a,b ; \b; \mu} \)=O\left(n^{-1/2}\right),\mbox{ as }n\rightarrow\infty.
	$$
\end{thm}

\begin{rem}
{\rm If we apply a partition $\g$ as introduced in the previous section directly, then the bound we may obtain is at most of order $o(1)$. }
\end{rem}

\noindent{\it Proof.}  According to Remark~\ref{meanchange}, with $\nu(dx)=\frac{1}{\ex|\Xi_1|}\ex\Xi_1(dx)$, it suffices to show
\begin{equation}
d_2 \( \L \( \sum_{k=1}^n \Xi_k \) , \bpi_{a,b ; \b; \nu} \)=O\left(n^{-1/2}\right),\mbox{ as }n\rightarrow\infty.\label{Bernoullishiftproof0}
\end{equation} 
We embed $\{\Xi_i\}$ into marked point processes \cite[pp.~194--195]{Daley08} and use Theorem~\ref{3.1} to complete the proof. To this end, 
we take fixed points $\{t_i:\ 1\le i\le m\}$ with $d_0'$ distances $1$ from each other and define a {\it ground process} \cite[p.~194]{Daley08} of $\Xi_1$ as $\Xi'_1=\sum_{i=1}^nI_i\d_{t_i}$. 
The metric $d_0'$ induces $d_1'$ and $d_2'$ in the same way as that $d_0$ generates $d_1$ and $d_2$. The mean measure of $\Xi'_1$ is $\lambda'=\sum_{i=1}^mp_i\d_{t_i}$. Let $\theta_1,\theta_2;a,b;\b$ be the same as those defined in Theorem~\ref{3.2} and set
$\nu'(dx) =\left(\sum_{i=1}^mp_i\right)^{-1}\lambda'(dx).$
For $\iid$ Bernoulli processes $\{\Xi'_k\}$, it follows from Theorem~\ref{3.1}, Remark~\ref{Bernoullinorole} and Remark~\ref{meanchange} that 
\begin{equation}
	d_2' \( \L \( \sum_{k=1}^n \Xi_k' \) , \bpi_{a,b ; \b; \nu'} \)=O\left(n^{-1/2}\right)\mbox{ as }n\to\infty.\label{Bernoullishiftproof1}
	\end{equation}
Using the Rubinstein duality theorem \cite[Theorem 8.1.1]{Rachev} and decompositions of point processes \cite[\S2.1]{Kallenberg83}, we can find $\Z_+^m$-valued random vectors $\(Y_{11},\dots,Y_{1m}\)$ and $\(Y_{21},\dots,Y_{2m}\)$ such that $\L\left(\sum_{i=1}^mY_{1i}\d_{t_i}\right)=\L\left(\sum_{k=1}^n \Xi_k'\right)$, 
$\L\left(\sum_{i=1}^mY_{2i}\d_{t_i}\right)=\bpi_{a,b ; \b; \nu'}$ and 
 \begin{equation}\ex d_1'\left(\sum_{i=1}^mY_{1i}\d_{t_i},\sum_{i=1}^mY_{2i}\d_{t_i}\right)=d_2' \( \L \( \sum_{k=1}^n \Xi_k' \) , \bpi_{a,b ; \b; \nu'} \).\label{Bernoullishiftproof2}\end{equation}
We now use $\sum_{i=1}^mY_{1i}\d_{t_i}$ and $\sum_{i=1}^mY_{2i}\d_{t_i}$ as ground processes
 to construct marked point processes as suitable realisations of $\L \( \sum_{k=1}^n \Xi_k \)$ and $ \bpi_{a,b ; \b; \nu}$. 
 Let $\{\(\zeta_{1i},\dots,\zeta_{mi}\):\ i\in \N\}$ be independent copies of $\(\zeta_1,\dots,\zeta_m\)$ such that $\{\(\zeta_{1i},\dots,\zeta_{mi}\):\ i\in \N\}$ is independent of
 $\{\(Y_{11},\dots,Y_{1m}\),\(Y_{21},\dots,Y_{2m}\)\}$, define
  $$\cw_1:=\sum_{i=1}^m\sum_{j=1}^{Y_{1i}}\d_{\zeta_{ij}}\mbox{ and }\cw_2:=\sum_{i=1}^m\sum_{j=1}^{Y_{2i}}\d_{\zeta_{ij}},$$
 then $\L\(\cw_1\)=\L \( \sum_{k=1}^n \Xi_k \)$, $\L\(\cw_2\)=\bpi_{a,b ; \b; \nu}$, and
 \begin{equation}
 d_2 \( \L \( \sum_{k=1}^n \Xi_k \) , \bpi_{a,b ; \b; \nu} \)\le \ex d_1\(\cw_1,\cw_2\)\le \ex d_1'\left(\sum_{i=1}^mY_{1i}\d_{t_i},\sum_{i=1}^mY_{2i}\d_{t_i}\right).\label{Bernoullishiftproof3}\end{equation}
 Combining \Ref{Bernoullishiftproof1}, \Ref{Bernoullishiftproof2} and \Ref{Bernoullishiftproof3} gives \Ref{Bernoullishiftproof0}.
 \qed

\subsection{Compound Poisson process}\label{sectcompoundPoisson}

\cite{BCL92} and \cite{BM02} demonstrate that a compound Poisson process is often good enough as a suitable asymptotic model for a variety of
random phenomena. In this example, we show that the superposition of such a model can be well described by Theorem~\ref{mainresult}.

Recall that a compound Poisson process on a compact carrier space $\G$ is defined as $\Xi_1 = \sum_{i=1}^\infty i X_i$, where $\{X_i\}$ are independent Poisson processes with mean measures $\{\lambda_i\}$ on $\G$ respectively and we write $\Xi_1\sim{\rm CP}(\lambda_1, \lambda_2,\dots)$.
\ignore{For brevity, we write $\Xi  = \sum_{i=1}^\infty i W_i \sim{\rm CP}(\lambda_1,\lambda_2,\dots)$, where $W_i$ is the sum of $n$ independent copies of $X_i$. Note that its superposition is still a compound Poisson process $\Xi\sim{\rm CP}(n \lambda_1, n \lambda_2,\dots)$, consequently, with suitably chosen parameters, $\bpi_{a,b;\b;\mu}$ can be used to replace a compound Poisson process in the context of superposition of point processes.}

	\begin{thm} \label{3.4} If $\Xi_1\sim{\rm CP}(\lambda_1, \lambda_2,\dots)$ with $\ex|\Xi_1|<\infty$ and $\sum_{j\ge 2}\lambda_j$ being absolutely continuous with respect to $\lambda_1$, then for \iid\ compound Poisson processes $\{\Xi_k\}_{k \in \mathbb{N}}$, with 
	$a,b;\b;\mu$ chosen as in Case 1, we have
	$$
	d_2 \( \L \( \sum_{k=1}^n \Xi_k \) , \bpi_{a,b ; \b; \mu} \)=o(1),\mbox{ as }n\to\infty.$$
\end{thm}

\begin{rem}
{\rm Noting that the superposition $\sum_{i=1}^n\Xi_i\sim{\rm CP}(n \lambda_1, n \lambda_2,\dots)$, Theorem~\ref{3.4} states that, with suitably chosen parameters, $\bpi_{a,b;\b;\mu}$ can be used to replace a compound Poisson process in the context of superposition of point processes.}
\end{rem}

\begin{rem}
{\rm The condition that $\sum_{j\ge 2}\lambda_j$ is absolutely continuous with respect to $\lambda_1$ guarantees aperiodicity of the distribution and it plays the crucial role 
in the theory of compound Poisson approximation in \cite{BCL92,BU98,BU99,BM02,Xia05a}.}
\end{rem}

\begin{rem}
{\rm It can be observed from the proof below that better upper bounds are possible if more information about $\{\lambda_i\}$ is available.}
\end{rem}

\noindent{\it Proof of Theorem~\ref{3.4}.}  Taking a reduced carrier space if necessary, without loss of generality, we assume $\G$ equals the support of $\lambda_1$, 
that is, the smallest closed set $A$ such that $\lambda_1(A)=|\lambda_1|$.  Since $\var(|\Xi_1|)=\sum_{i=1}^\infty i^2\lambda_i(\G)\ge\ex|\Xi_1|=\sum_{i=1}^\infty i\lambda_i(\G)$, Case~1 applies. Set $\Xi  = \sum_{i=1}^{n-2}\Xi_i$, let $\g = \{ G_1, \cdots, G_k \}$ be a partition and $W_1$ be a Poisson process on $\G$ with mean measure $(n-2)\lambda_1$, then
 \beas
& & d_{tv} (\L (\M_\g \circ \Xi + \d_{t_j}) , \L (\M_\g \circ \Xi )) \le d_{tv} (\L (\Xi (G_j) + 1), \L (\Xi (G_j)))
  \\ & \le &
d_{tv} (\L (W_1 (G_j) + 1), \L (W_1 (G_j)) ) \le \frac 1 {\sqrt{ 2e (n-2) \lambda_1 (G_j)}} ,
 \enas
 where the last inequality is from Proposition~A.2.7 in \cite{BHJ}.
Hence
 $$
\tv_\g (\Xi) \le \max_{1\le j\le k} \frac 1 {\sqrt{ 2e (n-2) \lambda_1 (G_j)}}, \ \ \ n \ge 4,
 $$
which implies, for arbitrary $\e>0$, $\vartheta_\e (\Xi) \le \inf_{\g\in\P_\e} {\max_{G\in\g} {\frac 1 {\sqrt{ 2e (n-2) \lambda_1 (G)}}}}$. It then follows from \Ref{mainresultineq1.1} that
$$\limsup_{n\to\infty}d_2 \( \L \( \sum_{k=1}^n \Xi_k \) , \bpi_{a,b ; \b; \mu} \)\le 2\e,$$
completing the proof. 
\qed

\subsection{Renewal process}\label{sectrenewal}

The superposition of renewal processes is not a renewal process except that they are Poisson processes \cite[p.~370]{Feller68} and the exact behaviour of the superposition is generally hard to extract. In this subsection, we establish its asymptotic behaviour.

Let $W_0$, $W_1$, $W_2$, $\cdots$ be independent non-negative random variables defined on a probability space $(\Omega, \F, \mathbb{P})$. The variables $W_1$, $W_2$, $\cdots$ are strictly positive and identically distributed, which play the role of {\it inter-renewal times} of the renewal process $\mathbf{S}=(S_n)_0^\infty:=(\sum_{i=0}^n W_i)_0^\infty$. We assume $\ex(W_1^2)<\infty$ and choose the delay $W_0$ to make the renewal process stationary \cite[p.~75]{Daley08}.  	
	We define $\Xi_1=\sum_{m=0}^\infty \d_{S_m} \mathbf{1}_{S_m\in \G}$, which is the 
	renewal point process {\it restricted} to $\Gamma=[0,1]$ \cite[p.~12]{Kallenberg83}. Before stating the result in this subsection, we briefly recall three terminologies. The support of a random variable $X$ is defined as the smallest closed set $A$ such that $\prob(X\in A)=1$ and, for two subsets $B_1,B_2$ of $\re$, $B_1+B_2:=\{x+y:\ x\in B_1,y\in B_2\}$ and $d_0(B_1,B_2)=\inf\{d_0(x,y):\ x\in B_1,y\in B_2\}$. 
	
	\begin{thm} \label{3.3} Assume the renewal time $W_1$ satisfies 
	\begin{equation}d_0( \supp(W_1)+\supp(W_1),\supp(W_1))=0\label{renewalcondition}
	\end{equation}
	and $\var(|\Xi_1|)> \frac12\ex|\Xi_1|$, then for $\iid$ renewal processes $\{\Xi_k\}_{k \in \mathbb{N}}$, with 
	$a,b;\b;\mu$ chosen as in Case 1 if $\var(|\Xi_1|)\ge \ex|\Xi_1|$ and in Case 2 if $\frac12\ex|\Xi_1|< \var(|\Xi_1|)<\ex|\Xi_1|$, we have
		$$
		d_2 \( \L \( \sum_{k=1}^n \Xi_k \) , \bpi_{a,b ; \b; \mu} \)=o(1), \mbox{ as }n\rightarrow\infty.
		$$
	\end{thm}

\begin{rem} {\rm If $0 \in \supp(W_1)$, then it satisfies \Ref{renewalcondition} and the bound in Theorem~\ref{3.3} holds.}
	\end{rem}
	
\begin{rem} {\rm The condition \Ref{renewalcondition} is almost necessary. See counterexample~\ref{counterex1} below.}
	\end{rem}
	
\begin{rem}{\rm The condition $\var(|\Xi_1|)> \frac12\ex|\Xi_1|$ can not be easily deduced from the moments of $W_1$. However, if we consider a sufficiently large carrier space, the asymptotic behaviour of the renewal process ensures that the condition can be verified through the first two moments of $W_1$.}
\end{rem}

\noindent{\it Proof of Theorem~\ref{3.3}.} For any $\epsilon>0$, we take an $m$ such that $2^{-m}<\epsilon$. We divide $\G$ into $2^m$ equally spaced intervals with $s_j = j / 2^m$ so that $\g = \{G_1,\dots,G_{2^m}\}$, where $G_1=[0,s_1]$ and $G_j = \( s_{j-1}, s_j \]$ for $2 \le j\le
2^m$. The centre of $G_j$ is $t_j = (s_{j-1} + s_j)/2$ and $d_0(\g) = 2^{-(m+1)}$. Consequently, the first term in \Ref{mainresultineq1.1} is bounded by $\epsilon$. With the partition $\g$, set $X_j = \Xi (G_j)$, define $\vec X = (X_1, \cdots , X_{2^m})$ and $\vec Y$ as the sum of $n-2$ independent copies of $\vec X$. Applying \cite[Lemma 4.1]{BLX17}, we obtain
  \beas
	d_{tv} (\L (\vec Y + e_j), \L (\vec Y ))  \le  \frac{C}{\sqrt{(n-2) u_m}},
\enas  where $u_m:=\min_{1\le j \le 2^m}\{1-d_{tv}(\L(\vec{X}), \L(\vec{X}+e_j) )\}$ and $C$ is a universal constant. If $u_m\neq0$, then
the second term in \Ref{mainresultineq1.1} with $\e=2^{-(m+1)}$ is also dominated by $\epsilon$ for sufficiently large $n$, which implies that the bound
in \Ref{mainresultineq1.1} can be made arbitrarily small {as $n\to\infty$}. To establish $u_m\neq0$,
we make use of the assumption that the support $A$ of $W_1$ satisfies $d_0(A+A,A)=0$. Since $A$ is closed, and in $\re_+$, the operation $+$ is continuous, $A+A$ and $A$ are both closed, which means that $(A+A)\cap A\neq \emptyset $. This in turn implies that there exists at least one $x \in \mathbb{R}_+$ such that both $\prob(W \in (x-\epsilon_1,x+\epsilon_1))$ and $\prob(W_1+W_2 \in (x-\epsilon_1,x+\epsilon_1))$ are positive for all $\epsilon_1 >0$. It is also possible to find a $0<y<x$ such that $\mathbb{P}(W_1 \in (y-\epsilon_2,y+\epsilon_2), W_1+W_2\in  (x-\epsilon_1,x+\epsilon_1))\neq 0$ for all $\epsilon_1$, $\epsilon_2>0$. For the convenience of argument, we extend the stationary renewal point process $\Xi_1$ to $\Xi_1'$ on $\re$. For $0< j\leq m$, if $\varsigma>0$ is small enough, the set
$$B_\varsigma:=\left(\left.\frac{j-1}{2^m}+\varsigma-y,\frac{j}{2^m}-\varsigma-y\right] \right\backslash
\left\{\bigcup_{n_1\leq 2^m}\bigcup_{n_2\in \mathbb{N}}\left(\frac{n_1}{2^m}-2\varsigma -n_2x,\frac{n_1}{2^m}+2\varsigma-n_2x\right]\right\}$$ has positive Lebesgue measure. 

\begin{figure}
\begin{center} 
 \begin{tikzpicture}[scale=1]
 \draw (0,0) -- (14,0);
 \draw (3,0) -- (3,0.2);
 \draw (5,0) -- (5,0.2);
 \draw (7,0) -- (7,0.2);
 \draw (9,0) -- (9,0.2);
 \draw (11,0) -- (11,0.2);
 \node[] at (3,0.35)   () {0};
 \node[] at (5,0.35)   () {0.25};
 \node[] at (7,0.35)   () {0.5};
 \node[] at (9,0.35)   () {0.75};
 \node[] at (11,0.35)   () {1};
 \node[circle,draw,inner sep=2pt] at (2,0.6)   (a1) {};
 \node[circle,draw,inner sep=2pt] at (6,0.6)   (a2) {};
 \node[circle,draw,inner sep=2pt] at (10,0.6)   (a3) {};
 \node[circle,draw,inner sep=2pt] at (14,0.6)   (a4) {};
 \draw [-](a1)--node[above]{$\approx0.5$} (a2);
 \draw [-](a2)--node[above]{$\approx0.5$} (a3);
 \draw [-](a3)--node[above]{$\approx0.5$} (a4);
 \node[regular polygon,regular polygon sides=3,draw,inner sep=1.5pt] at (2,-0.2)   (b1) {};
 \node[regular polygon,regular polygon sides=3,draw,inner sep=1.5pt] at (3.6,-0.2)   (b2) {};
 \node[regular polygon,regular polygon sides=3,draw,inner sep=1.5pt] at (6,-0.2)   (b3) {};
 \node[regular polygon,regular polygon sides=3,draw,inner sep=1.5pt] at (10,-0.2)   (b4) {};
 \node[regular polygon,regular polygon sides=3,draw,inner sep=1.5pt] at (14,-0.2)   (b5) {};
 \draw [-](b1)--node[below]{$\approx0.2$} (b2);
 \draw [-](b2)--node[below]{$\approx0.3$} (b3);
 \draw [-](b3)--node[below]{$\approx0.5$} (b4);
 \draw [-](b4)--node[below]{$\approx0.5$} (b5);
 \end{tikzpicture}
\begin{minipage}{0.70\textwidth}\centering
\caption{\label{fig1} $m=2$, $x=0.5$, $y=0.2$, $j=1$}
\end{minipage}
\end{center}
\end{figure}
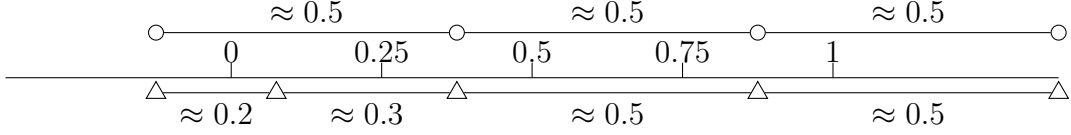

From stationarity, there is a positive probability that there is at least one point in $B_\varsigma$, and conditional on the largest point in $B_\varsigma$ and the past, the renewal process has a positive probability for the future inter-renewal times $W_1'$, $W_2'$, $\cdots$ to evolve as $W_{i}' \in (x-\frac{2\varsigma}{2^i},x+\frac{2\varsigma}{2^i})$ for all $i\in \mathbb{N}$ until time $1$ and it also guarantees a positive probability that the incoming inter-renewal times $W_1''$, $W_2''$, $\cdots$ evolve as $W_1'' \in (y-\varsigma,y+\varsigma)$, $W_1''+W_2''\in  (x-\varsigma,x+\varsigma)$, ${W_i'' =W_{i-1}'}\in (x-\frac{4\varsigma}{2^i},x+\frac{4\varsigma}{2^i})$ for $i\ge 3$ until time $1$. The choice of $B_\varsigma$ and synchronicity of $\{W_i':\ i\ge 1\}$ and $\{W_i'':\ i\ge 1\}$ ensure that an extra renewal point caused by $W''_1$ is added in $G_j$, and the subsequent renewal points of the two renewal processes occur in the same partition sets $\{G_k:\ k\ge j\}$ simultaneously.
{Consequently,} we can set aside a positive probability event $B_+$ such that on $B_+$, {the} two renewal processes run together until the point in $B_\varsigma$ and then one runs according to $\{W_i'\}$ and the other evolves as $\{W''_i\}$. 
Figure~\ref{fig1} shows the coupling when $m=2$, $x=0.5$, $y=0.2$, $j=1$,
a renewal happens at around $-0.25$, with a positive probability, the next three inter-arrival times are each around $x=0.5$; with another positive probability, the incoming four inter-arrival times respectively take values around $y=0.2$, $x-y=0.3$, $x=0.5$, $x=0.5$. 
For this coupling, the corresponding vectors $\vec X'$ and $\vec X''$ satisfy that $\vec X''=\vec X'+e_j$ on $B_+$ (in Figure~\ref{fig1}, $\vec X'=(0,1,0,1)$ and $\vec X''=(1,1,0,1)$),
which implies that $d_{tv}(\L(\vec{X}), \L(\vec{X}+e_j) )<1$ for all $0<j\le 2^m$. This concludes the proof. \qed   

\begin{counterex}\label{counterex1}{\rm If $\supp{(W_1)}\subset[a,2a-b]$ for some $0<b<a\le \frac{1}{2}$, then $u_m=0$ for some $m$ so the method does not work.}
	\end{counterex}
	
In fact, for $m$ large enough, when $X_{2^{m-1}}=0$, there is one point of $\mathbf{S}$ sitting in the interval $B_3:=\left[\frac{1}{2}-\frac{1}{2^m}-a+\frac{b}{2},\frac{1}{2}+a-\frac{b}{2}\right]$ almost surely. But when we have $X_{2^{m-1}}=1$, there are no points in $B_3$ except in $G_{2^{m-1}}$. On the other hand, for $m$ large enough, $X_{2^{m-1}}\leq 1$ almost surely because $a>0$. In this situation, $d_{tv}(\L(\vec{X}), \L(\vec{X}+e_{2^{m-1}}) )=1$, i.e., $u_m=0$.

\begin{rem}{\rm It is possible to extend Theorem~\ref{3.3}
to the superposition of $\iid$ non-stationary renewal processes, provided there are use-friendly criteria for ensuring $u_m:= \min_{1\le j \le 2^m}\{1-d_{tv}(\L(\vec{X}), \L(\vec{X}+e_j) )\}\neq0$ for all $m\ge 1$ and $0<j\le 2^m$.}
\end{rem}

\subsection{Entrances or exits of Markov Process }\label{sectmarkov}

Let $\{M_t\}_{t\in\mathbb{R}}$ be a
time-reversible and irreducible Markov chain with finite state space $E$. Let $S_0$ be a proper subset of $E$. As the exit process from $S_0$ can be viewed as the entrance process of $E\backslash S_0$, we  consider entrance process to $S_0$ only. Let $T_1:=\inf \{t\ge 0:\ M_{t^-} \notin S_0 \mbox{ and }M_t \in S_0\}$ and $T_{i+1}=\inf \{t> T_{i}: M_{t^-} \notin S_0 \mbox{ and }M_t \in S_0\}$ for $i\geq 1$. Then the total number of entrances to $S_0$ in $\G=[0,1]$ can be written as $\tau=\max\{n: T_n\leq 1\}$ with $\max\emptyset:=0$, and the times of entrances form a point process $\Xi_1:=\sum_{1\leq i\leq \tau} \d_{T_i}$ with convention $\Xi_1=0$ when $\tau=0$. 
Clearly, $|\Xi_1|$ is almost surely finite. 

	\begin{thm} \label{entrance} For \iid\ entrance processes $\{\Xi_k\}_{k \in\N}$, with 
	$a,b;\b;\mu$ chosen as in Case 1,
		$$
		d_2 \( \L \( \sum_{k=1}^n \Xi_k \) , \bpi_{a,b ; \b; \mu} \)=o(1), \mbox{ as }n\to\infty.		$$
	\end{thm}

\begin{rem}{\rm When $S_0$ is a single point set, $(T_n)_{0}^\infty$ forms a renewal process, Theorem~\ref{entrance} becomes a special case of Theorem~\ref{3.3}.	However, when $S_0$ contains more than one state, then $\Xi_1$ is no longer a renewal process.}
\end{rem}

\noindent{\it Proof of Theorem~\ref{entrance}.}  \cite[Corollary~2]{BHX98} implies that $\var(|\Xi_1|)\ge \ex|\Xi|$, so Case 1 applies.
The rest of the proof is essentially the same as that of Theorem~\ref{3.3}. For any $\epsilon>0$, we choose an $m$ such that $2^{-m}<\epsilon$. Let $s_j = j / 2^m$ and $\g = \{G_1,\dots,G_{2^m}\}$, where $G_1=[0,s_1]$ and $G_j = \( s_{j-1}, s_j \]$ for $2 \le j\le
2^m$. The centre of $G_j$ is at $t_j = (s_{j-1} + s_j)/2$ and $d_0(\g) = 2^{-(m+1)}$. This partition ensures that the first term in \Ref{mainresultineq1.1} is bounded by $\epsilon$. Set $X_j = \Xi (G_j)$ and define $\vec X = (X_1, \cdots , X_{2^m})$ and $\vec Y$ as the sum of $n-2$ independent copies of $\vec X$. It follows from \cite[Lemma 4.1]{BLX17} that
  \beas
	d_{tv} (\L (\vec Y + e_j), \L (\vec Y ))  \le  \frac{C}{\sqrt{(n-2) u_m}},
\enas  where $u_m:=\min_{1\le j \le 2^m}\{1-d_{tv}(\L(\vec{X}), \L(\vec{X}+e_j) )\}$ and $C$ is a universal constant. It remains to show 
that $u_m\neq0$. Since $\{M_t\}_{t\in\re}$ is irreducible, we can choose a state $s\in E\backslash S_0$ such that there is a positive probability of entering $S_0$ {immediately after leaving} $s$. Let $\tau_1$ be the first time that the Markov chain enters $S_0$, $\tau_2$ be the first time after $\tau_1$ to depart from $S_0$, $T_1'$ and $T_2'$ be the first and second jump times of $\{M_t\}_{t\in\mathbb{R}}$ after time $0$. From the assumption that $\{M_t\}_{t\in\re}$ is finite irreducible, we can conclude that $\prob(M_0=s)>0$ and so $q_1:=\prob\(\vec X= \mathbf{0}\)\ge \prob\(M_0=s,\tau_1>1\)\ge\prob\(M_0=s,T_1'>1\)>0$ and 
\begin{eqnarray*}
q_2&:=&\prob\(\vec X=e_j\)\ge\prob\(M_0=s,\tau_1\in G_j,\tau_2>1\)  \\
&\ge&\prob\(M_0=s,T_1'\in G_j, M_{T_1'}\in S_0,T_2'>1\)>0,\end{eqnarray*}
which in turn imply 
 $d_{tv}(\L(\vec{X}), \L(\vec{X}+e_j) )\le \max\{1-q_1,1-q_2\}<1$, 
 as claimed. \qed

\section{The Proof of Theorem~\ref{mainresult}} \label{Section.Proof}

The advantage of using $\bpi_{a,b;\b;\mu}$ as approximating distribution is that it can be considered as 
the unique stationary distribution of an $\H$-valued positive recurrent process with the generator
 \begin{eqnarray*}
\A h (\xi)&:=& \big( a + b |\xi| \big)  \int_\G \big( h (\xi + \d_x) - h (\xi)
\big) \mu (dx) \\
&&+ \big( 1 + \b( |\xi| - 1) \big) \int_\G \big( h(\xi - \d_x) -
h(\xi) \big) \xi (d x),
 \end{eqnarray*}
 see \cite{XZ12} for more details.
We use $\Z_\xi(\cdot) $ to stand for a birth-death point process with generator $\A$ and initial configuration $\xi$. For any bounded measurable function $f$ on $(\H,\B(\H))$, it can be shown that 
$$h_f (\xi) := - \int_0^\infty \big( \ex f (\Z_\xi (t)) - \bpi_{a,b;\b;\mu} (f)
\big) d t$$ 
is well defined and is the solution of the Stein equation
 $$
\A h(\xi) = f(\xi) - \bpi_{a,b;\b;\mu} (f) .
 $$
To estimate $d_2 (\L (\cw), \bpi_{a,b;\b;\mu})$, it is equivalent to bound $\ex \A h_f (\cw)$ for all $f \in \F$
defined on page~\pageref{pageofcapitalF}. As $\ex \A h (\cw)$ can be expressed via the differences of $h$, the successful application of the Stein method hinges on sharp upper bounds of
 $$
\D h (\xi; x) : =  h (\xi + \d_x) - h (\xi) , \ \D_2 h (\xi; x,y) : = \D h (\xi + \d_y; x) - \D h (\xi; x) .
 $$
Let $\D_2 h (\xi) := \sup \{ |\D_2 h(\xi; x, y)| :\ x,y \in \G \}$. Then it is shown in \cite{XZ12} that
 \eq \label{Eq.Steinfactor}
\D_2 h_f (\xi) \le \frac 2 {| \xi | + 1} + \frac 5 a , \ \ \ \forall f \in \F, \xi \in \H.
 \en

Now, we are ready to prove Theorem~\ref{mainresult}.

\noindent{\it Proof of Theorem~\ref{mainresult}.} The inequalities \Ref{mainresultineq1.1} and \Ref{mainresultineq2.1} are due to the well-known concentration inequality, see \cite[Theorem~2.7]{Mc98}
and \cite[Theorem~2.7]{CL06}. Hence it remains to show \Ref{mainresultineq1} and \Ref{mainresultineq2}.

Suppose $\g \in \P_\e$. The ``assembling mapping" $\M_\g$ ensures that for any configuration $\eta$,
 $$
d_1 (\eta, \M_\g \circ \eta) \le d_0 (\g) \le \e .
 $$
It follows that $ d_2 (\L (\cw), \L (\M_\g \circ \cw)) \le \e$ for any point process $\cw$, which yields
 \bea
 & &
d_2 (\L (\Xi), \bpi_{a,b;\b;\mu}) \nonumber
 \\ & \le &
d_2 (\L (\Xi),  \L (\M_\g \circ \Xi)) + d_2 (\L (\M_\g \circ \Xi), \L (\M_\g \circ \cz)) \nonumber\\
&&+ d_2 (\L (\M_\g \circ \cz), \bpi_{a,b;\b;\mu}) 
 \\  & \le &
2 \e + d_2 (\L (\M_\g \circ \Xi), \L (\M_\g \circ \cz)) . \label{Eq.proof1}
 \ena

To compute the $d_2$ distance between $\L (\M_\g \circ \Xi)$ and $\L (\M_\g \circ \cz)$, we concentrate on the space $\tilde \G := \{ t_1, \cdots, t_k \}$ and apply Stein's method. Denote by $\tilde \H$ the class of all configurations on $\tilde \G$.
For any $\xi \in \tilde \H$, let
 $$
\tilde \A \tilde h (\xi) : = (a + b |\xi|) \int \D \tilde h (\xi; x) \tilde \mu (d x) - (1 + \b (|\xi|-1)) \int \D \tilde h (\xi - \d_x; x) \xi (d x) .
 $$
Then, with generator $\tilde \A$, we have a positive recurrent Markov process on $\tilde \H$. The unique stationary measure is $\tilde \bpi : = \bpi_{a,b;\b; \tilde \mu} = \L (\M_\g \circ \cz)$, where
 $$
\tilde \mu (dx)= \sum_{i=1}^k \mu (G_i) \d_{t_i} (dx) .
 $$
Denote
 $$
\tilde \F : = \{ \tilde f : |\tilde f (\xi) - \tilde f (\eta) | \le d_1 (\xi, \eta), \forall \xi, \eta \in \tilde \H \} .
 $$
For any $\tilde f \in \tilde \F$, let $\tilde h_{\tilde f}$ be the unique solution of
 $$
\tilde \A \tilde h_{\tilde f} = \tilde f - \tilde \bpi ( \tilde f) .
 $$
Then
 \eq \label{Eq.proof2}
d_2  (\L (\M_\g \circ \Xi), \L (\M_\g \circ \cz)) = \sup_{\tilde f \in \tilde \F} \ex \tilde \A \tilde h_f (\M_\g \circ \Xi) .
 \en

Now we concentrate on estimating $\ex \tilde \A \tilde h_f (\M_\g \circ \Xi)$. Let
 $$
h (\xi) : = \tilde h_{\tilde f} (\M_\g \circ \xi) , \ \ \ \forall \xi \in \H.
 $$
Then
 $$
\tilde \A \tilde h (\M_\g \circ \xi) = (a + b |\xi|) \int \D h (\xi; x) \mu (d x) - (1 + \b (|\xi| - 1)) \int \D h (\xi - \d_x; x) \xi (d x) .
 $$

First of all, we can write $\ex \tilde \A \tilde h (\M_\g \circ \Xi)$ via $\D^2 h$'s. Namely, since $\Xi_1, \cdots, \Xi_n$ are independent identically distributed,
 \beas
&&\ex \tilde \A \tilde h (\M_\g \circ \Xi)\\
 &&= a \int \ex \[ \D h (\Xi; x) \] \mu (dx) + b \int \ex \[ \Xi \D h (\Xi; x) \] \mu (dx)
 \\ & &\ \ \   -  \ex \[ \( 1 + \b (|\Xi| - 1) \) \int \D h (\Xi - \d_x; x) \Xi (dx) \]
\\ &&= a \int \ex \[ \D h (\Xi; x) \] \mu (dx) + b n \int \ex \[ \Xi_1 \D h (\Xi; x) \] \mu (dx)
 \\ &&\ \ \  - n \ex \[ \(1 +  \b (|\Xi_1| - 1) + \b (n-1) |\Xi_2| \) \int \D h (\Xi - \d_x; x) \Xi_1 (dx) \] .
 \enas
With the reduced Palm processes, one can write $\ex \tilde \A \tilde h (\M_\g \circ \Xi)$ as
 \bea
 && a \int \ex \[ \D h (\Xi; x) \] \mu (dx) + b n \iint \ex \[ \D h \left(\Xi_{1,y} + \d_y + \snn ; x\right) \] \mu (dx) \l (dy) \nonumber
 \\ & & - n \int  \ex \[ \D h \left(\Xi_{1,x} + \snn ; x\right) \] \l (dx) \nonumber
 \\ & & - \b n \iint \ex \[ \D h \left(\Xi_{1,x,y} + \d_y + \snn; x\right) \] \l^{[2]} (dx,d y) \nonumber
 \\ & & - \b n (n-1) \iint \ex \[ \D h \left(\Xi_{1,x} + \Xi_{2,y} + \d_y + \snnn ; x\right) \] \l (dx) \l (dy) . \label{Eq.intoDh}
 \ena
We subtract $\D h (\snn; x)$ in the first four terms and $\D h (\snnn ; x)$ in the last one. Then,  $\ex \tilde \A \tilde h (\M_\g \circ \Xi)$ can be written via $\D_2 h$'s, provided that the number of $-\D_2 h$ added is balanced with that of $\D_2$ added. More precisely, we need
 $$
a \mu (dx) + b n |\l| \mu (dx) - n \l (dx)  - \b n \l^{[2]} (dx, \G) = \b n (n-1) |\l| \l (dx) ,
 $$
which is equivalent to \eqref{Eq.a} and \eqref{Eq.mu}. With \eqref{Eq.a} and \eqref{Eq.mu}, we write $\ex \tilde \A \tilde h (\M_\g \circ \Xi)$ via $\D_2 h$'s. For example, the first term in \eqref{Eq.intoDh} becomes
 $$
a \int \ex \[ \D h (\Xi; x) - \D h ( \snn ; x) \] \mu (dx) .
 $$
The difference $\D h (\Xi; x) - \D h ( \snn ; x)$ can be telescoped out as the sum of $|\Xi_1|$ $\D_2 h$ functions. Provided the number of $\D_2 h$ is balanced with that of $- \D_2 h$, one can further write $\ex \tilde \A \tilde h (\M_\g \circ \Xi)$ via $\D_3 h$'s or differences of two $\D_2 h$'s. To this end, let $z \in \G$ and

 	\beas
e_1(x) & = & \ex  \[ \D h (\Xi; x) - \D h \left(\snn; x\right)  \] - \ex \[ |\Xi_1| \D_2 h \left(\snn; z,z\right)  \] ,
 \\
e_2(x,y) & = & \ex \[ \D h \left(\Xi_{1,y} + \d_y + \snn ; x\right) - \D h \left(\snn; x\right)  \]  \\
&&- \ex \[ (|\Xi_{1,y}| + 1 ) \D_2 h \left(\snn ; z,z\right) \] ,
 \\
e_3(x) & = & \ex \[ \D h \left(\Xi_{1,x} + \snn ; x\right) - \D h \left(\snn; x\right)  \] \\
&& - \ex \[ |\Xi_{1,x}|  \D_2 h \left(\snn ; z,z\right)  \] ,
 \\
e_4(x,y) & = & \ex \[ \D h \left(\Xi_{1,x,y} + \d_y + \snn; x\right) - \D h \left(\snn; x\right) \] \\
&&- \ex \[ ( |\Xi_{1,x,y}|+1) \D_2 h \left(\snn ; z,z\right)  \] ,
 \\
e_5(x,y) & = & \ex \[ \D h \left(\Xi_{1,x} + \Xi_{2,y} + \d_y + \snnn ; x\right) - \D h \left(\snnn ; x\right) \] \\
&& - \ex \[  (|\Xi_{1,x}| + |\Xi_{2,y}|+1)  \D_2 h \left(\snnn  ; z,z\right) \] ,
 \\
e_6(x) & = & \ex \[ \D h \left(\snn; x\right)  - \D h \left( \snnn ; x\right)  \] - \ex \[ |\Xi_2|   \D_2 h \left(\snnn  ; z,z\right) \] ,
 \\
e_7 & = & \ex \D_2 h \left(\snn; z,z\right) - \ex \D_2 h \left(\snnn ; z,z\right) .
\enas

Then,
 \beas
 &&
\ex \tilde \A \tilde h (\M_\g \circ \Xi)
 \\ &=& a \int e_1(x) \mu (dx)  + b n \iint e_2(x,y) \mu (dx) \l (dy) \\
 &&- n \int  e_3(x) \l (dx)  - \b n \iint e_4(x,y) \l^{[2]} (dx,d y)
 \\ & & - \b n (n-1) \iint e_5(x,y) \l (dx) \l (dy)  + \b n (n-1) |\l| \int e_6(x) \l (dx)
 \\ & &
+ \b n (n-1) ( 2 \theta_1 \theta_2 + \theta_1^2 - \theta_1^3 ) e_7 ,
 \enas
provided that
 \beas
 & &
 \theta_1^2  - \theta_2  + b \(  \theta_2 + \theta_1 - \theta_1^2 \)
 =
 \b (n-1)\theta_1 \[  \theta_1 - 2 ( \theta_1^2- \theta_2) \] + \b (\theta_3 + \theta_2 - \theta_1 \theta_2) .
 \enas
In both cases, $b$ and $\b$ are taken to ensure the above equality.

To estimate $e_1, \cdots, e_7$, we decompose them into the sum of $\D_2 h$ functions of the forms
 \beas
\D_3 h (x,y;z) & : = & \D_2 h (\xi + \d_z ; x,y) - \D_2 h (\xi; x,y), \\
D_{2,T} h (\xi; x,y; z,w) & := & \D_2 h (\xi; x,y) - \D_2 h (\xi; z,w) .
 \enas
The bounds in the following lemma can be found in \cite[pp.~3060-3061]{XZ12}.
 \begin{lem}
\label{Lem.D3}
For any point process $\cw$, and $u > 0$, both $|\ex \D_3 h (\cw; x,y,z) |$ and $ |\ex \D_{2,T}   h (\cw; x,y; z,z) | $ are bounded above by
 $$
r (\cw) := \frac {4u + 10}{a} \tv_\g (\cw) + 4 \prob \( 1 + |\cw| \le \frac a u \) ,
 $$
where $a$ is defined in \eqref{Eq.a} and $\tv_\g (\cw)$ is defined in \eqref{Eq.tv}.
 \end{lem}
To estimate $e_1$, let $\Xi_1 = \sum_{n=1}^{|\Xi_1|} \d_{X_n}$, $\langle \Xi_1\rangle_r : = \sum_{n=1}^r \d_{X_n}$, $\cw = \snn$. Then,
 \beas
 e_1(x) & = &
\ex  \[ \D h (\cw + \Xi_1 ; x) - \D h (\cw; x)  \] - \ex \[ |\Xi_1| \D_2 h (\cw; z,z)  \]
 \\ & = &
\ex \sum_{r=1}^{|\Xi_1|}  \[ \D_2 h (\cw + \langle \Xi_1\rangle_{r-1} ; x, X_r ) - \D_2  h (\cw; z,z)  \]
 \\ & = &
\ex  \sum_{r=1}^{|\Xi_1|}  \sum_{s=1}^{r-1}  \D^3 h (\cw + \langle \Xi_1\rangle_{s-1} ; x, X_r; X_s ) + \ex \sum_{r=1}^{|\Xi_1|}  D_{2,T}  h (\cw ; x, X_r; z,z)  .
 \enas
Since $\cw$ is independent of $\Xi_1$, it follows that
 $$
|e_1(x)| \le \( \frac 1 2 \ex \[ |\Xi_1| (|\Xi_1| - 1)  \] + \ex |\Xi_1| \) r (\cw) = \frac 1 2 \ex \[ |\Xi_1| (|\Xi_1| + 1)  \]  r \left(\snn\right).
 $$
Similarly, we have
 \beas
|e_2(x,y) | & \le &   \frac 1 2 \ex \[ (|\Xi_{1,y}| + 1)( |\Xi_{1,y}| +2) \]  r \left(\snn\right) , \\
|e_3(x)| & \le &  \frac 1 2 \ex \[  |\Xi_{1,x}| ( |\Xi_{1,x}| + 1) \] r \left(\snn\right)  , \\
|e_4(x,y)| & \le & \frac 1 2 \ex \[ (|\Xi_{1,x,y}| + 1) (|\Xi_{1,x,y}| + 2) \]  r \left(\snn\right) , \\
|e_5(x,y)| & \le & \frac 1 2 \ex \[ (|\Xi_{1,x}| + |\Xi_{2,y}| + 1) (|\Xi_{1,x}| + |\Xi_{2,y}| + 2) \]  r \left(\snnn \right) , \\
|e_6(x)| & \le & \frac 1 2 \ex \[ |\Xi_2| (|\Xi_2| + 1) \]  r \left(\snnn \right), \\
|e_7| & \le & \ex |\Xi_2| r \left(\snnn\right) .
  \enas
Since $r (\snn) \le r (\snnn)$, we have $ |\ex \tilde \A \tilde h (\Xi) |  \le C_n r (\snnn) $, where
 \beas
C_n & = & \frac 1 2 a ( \theta_2 + 2 \theta_1 ) + \frac 1 2 b n ( \theta_3 + 4\theta_2+2\theta_1 ) + \frac 1 2 n   ( \theta_3 + 2 \theta_2 ) + \frac 1 2 \b n ( \theta_4 + 4 \theta_3 + \theta_2 )
 \\ & & + \b n (n-1)  \left\{ \frac 1 2 \cdot 2\(\theta_1 (\theta_3 + 3 \theta_2 + \theta_1) + \theta_2 (\theta_2 + \theta_1) \) \right.\\
 &&\ \ \ \ \ \ \ \ \ \ \ \ \ \ \ \ \ \ \left.+ \frac 1 2 \theta_1^2 ( \theta_2 + 2 \theta_1  )  + \theta_1 |2 \theta_1 \theta_2 + \theta_1^2 - \theta_1^3| \right\} .
 \enas
It is not difficult to check that in each of the two cases, $a$ has order $n$, $b$ is a constant and $\b$ has order $1/n$. Hence $C_n$ has order $n$. Let $u$ be a constant independent of $n$ such that $ a / u < n \theta_1 / 4$, then for $n>2$,
 $$
|\ex \tilde \A \tilde h (\M_\g \circ \Xi)| \le C \cdot \tv_\g \left(\snnn\right) + Cn \prob \( \left|\snnn\right| \le \frac 3 4 \theta_1 (n-2) \) ,
 $$
where $C$ is a constant. Using the fact that $\snnn$ and $\sum_{k=1}^{n-2} \Xi_k$ have the same distribution, we combine \eqref{Eq.proof1} and \eqref{Eq.proof2} to conclude that
 \begin{eqnarray*}
&&d_2 \( \L \( \sum_{k=1}^n \Xi_k \) , \bpi_{a,b ; \b; \mu} \)\\
&& \le 2 \e + C\left\{ \tv_\g \( \sum_{k=1}^{n-2} \Xi_k \)  + n \prob \( \sum_{k=1}^{n-2} | \Xi_k | \le \frac 3 4 \theta_1 (n-2) \)\right\} .
 \end{eqnarray*}
Since $\g$ is arbitrary, the proof of Theorem~\ref{mainresult} is complete. \qed

 \end{document}